\documentclass[10pt]{article}

\usepackage{amssymb,amsmath,theorem}

\usepackage{color}                    
\definecolor{red}{rgb}{1,0,.2}        
\definecolor{cjp}{rgb}{.1,.7,.2}        

\newcommand{\N}{{\mathbb{N}}}

\advance \topmargin by -\headheight
\advance \topmargin by -\headsep
\evensidemargin \oddsidemargin
\author{
F. M. Dekking \\
Delft University of Technology \\
Faculty EEMCS, P.O.~Box 5031\\
2600 GA Delft, The Netherlands\\
{\tt F.M.Dekking@math.tudelft.nl}
}

\title{\bf Permutations of $\mathbb{N}$ generated by left-right filling algorithms}

\newcommand{\pq}[2]{ \left( \begin{smallmatrix} #1\\ #2 \end{smallmatrix} \right)}
\newcommand{\PQ}[2]{ \left( \begin{matrix} #1\\ #2 \end{matrix} \right)}
\newcommand{\proof} {\noindent {\it Proof:\;}}

 \theoremstyle{plain}
 \newtheorem{theorem}{Theorem}
 
 \newtheorem{lemma}[theorem]{Lemma}
 \newtheorem{proposition}[theorem]{Proposition}
 \theoremstyle{definition}

\begin{document}

\maketitle

\begin{abstract}
 We give an in depth analysis of an algorithm that generates permutations of the natural numbers introduced by Clark Kimberling in the On Line Encyclopedia of Integer Sequences. It turns out that the examples of such permutations in OEIS are completely determined by 3-automatic sequences.
\end{abstract}

\medskip


\section{Introduction}

Let $L:\N\rightarrow\N$ and $R:\N\rightarrow\N$ be two functions. A {\it left-right filling procedure} is an algorithm that may produce an permutation $\Pi$ of $\N$ as follows.\\
First of all $\Pi(1)=1$. Then for each $n\ge 2$, $$\Pi(n-L(n))=n\; {\;\rm if}\; \Pi(n-L(n)) \;{\rm is\; not\; yet\; defined,\; otherwise}\; \Pi(n+R(n))=n.$$
Here we say {\it may}, because $\Pi(n)$ might not be defined for all $n$.

\medskip

The simplest example of a pair $(L,R)$ generating a permutation is given by
$$L(n)=1, \;R(n) = 1 \quad {\rm for\; all\;} n.$$
We find consecutively: $n=2 \Rightarrow\Pi(2+R(2))=\Pi(3)=2$, $n=3 \Rightarrow\Pi(3-L(2))=\Pi(2)=3$, $n=4 \Rightarrow\Pi(4+R(4))=\Pi(5)=4$,
and in general, $\Pi(2n)=2n+1,\, \Pi(2n+1)=2n$. So this pair $(L,R)$ generates the self-inverse permutation of $\mathbb{N}$ given by sequence A065190 in OEIS (\cite{oeis}).

\medskip

An example of a pair $(L,R)$ that does not generate a permutation is given by
$$L(n)=n-1, \;R(n) = n-1 \quad {\rm for\; all\;} n.$$
Here only the odd numbers obtain a value: $\Pi(2n+1)=n+1$ for all $n\ge 0$.

We will be occupied mainly with the pair of functions given by
$$L(n)=\Big\lfloor \frac{n}2\Big\rfloor, \;R(n) = \Big\lfloor \frac{n}2\Big\rfloor \quad {\rm for\; all\;} n.$$
This pair was introduced by Clark Kimberling in the On Line Encyclopedia of Integer Sequences as A026136 (\cite{oeis}). The permutation it generates is
$$\Pi = 1, 3, 2, 7, 9, 4, 5, 15, 6, 19, 21, 8, 25, 27, 10, 11, 33, 12, 13, 39, 14, 43, 45, 16, 17, 51, 18, \dots$$
We will establish a one-to-one connection of this permutation with a 3-automatic sequence that permits to prove a number of properties of this $\Pi$, and related permutations.

\section{Left and right functions}\label{sec:AA}

\subsection{The pair $L(n)=\lfloor n/2\rfloor, \;R(n) = \lfloor n/2\rfloor.$}\label{subs:standard}

When $L(n)=\lfloor n/2\rfloor, \;R(n) = \lfloor n/2\rfloor$, then there are four possibilities for the value of $\Pi(n)$.
Splitting into the cases $L(n)=R(n)=n/2$ for $n$ even, and  $L(n)=R(n)=(n-1)/2$ for $n$ odd, one arrives at the following four cases\\[-.6cm] \begin{align*}
{\rm (I)}  &\quad \Pi(n)\; {\rm odd},\;\,  \; \Pi(n)>n,  \; \Pi(n) = 2n-1\\
{\rm (II)} &\quad \Pi(n)\; {\rm even},\,   \; \Pi(n)>n,  \; \Pi(n) = 2n\\
{\rm (III)} &\quad \Pi(n)\; {\rm odd},\;\, \; \Pi(n)<n,  \; \Pi(n) = (2n+1)/3\\
{\rm (IV)} &\quad \Pi(n)\; {\rm even},\,   \; \Pi(n)<n,  \; \Pi(n) = 2n/3.
\end{align*}

\noindent We say $\Pi(n)$ is of type {\rm (I)}, etc., when we are in case {\rm (I)}, etc.

\begin{lemma}\label{lem:step2n} If the left-right filling procedure is at step $2n$, then the natural numbers $1,2,\dots,n$ have already been assigned a $\Pi$-value.
\end{lemma}

 \proof \; By induction. This is trivially true for $n=1$. If we are at step $2(n+1)$, then $1,2,\dots,n$ have  been assigned a $\Pi$-value by the induction hypothesis.  The number $n+1$ might have got a $\Pi$-value in steps $1,2,\dots,2n$. If not, then it got the value $\Pi(n+1)=2n+1$ at step $2n+1$, since\\[-.2cm]
$$\hspace*{5cm}2n+1-\Big\lfloor\frac{2n+1}2\Big\rfloor=2n+1-n=n+1.\hspace*{5cm} \Box$$

 \medskip

\noindent Lemma \ref{lem:step2n} excludes the possibility to go to the left at step $2n$, so it implies the following proposition.

\begin{proposition}\label{prop:step2n} For all $n \ge 1$,  $\Pi(3n)=2n$.
\end{proposition}

\noindent The next task is to determine the type of all entries of $\Pi$.

\begin{theorem}  \label{th:types}
Let $n,k$ be  natural numbers. Then\\
{\rm \bf a)}\; $\Pi(n)$ never has type {\rm (II)}.\\
{\rm \bf b)}\; If\, $n=3k$ then $\Pi(n)$ has type {\rm (IV)}. \\
{\rm \bf c)}\; If\, $n=3k+2$ then $\Pi(n)$ has type {\rm (I)}.\\
{\rm \bf d)}\; If\, $n=3k+1$ then $\Pi(n)$ has type {\rm (I)} if $\Pi(k+1)$ has type {\rm (I)} and $\Pi(n)$ has type {\rm (III)} if $\Pi(k+1)$ has type {\rm (III)} or  type {\rm (IV)}.
\end{theorem}

\proof Part {\rm \bf a)}: $\Pi(n)$ is never equal to $2n$, since either $n$ has been assigned the value $2n-1$ at step $2n-1$:
$2n-1-\lfloor (2n-1)/2\rfloor=n$, or $n$ has been assigned a smaller value at an earlier step.

\noindent Part {\rm \bf b)} is a rephrasing of Proposition \ref{prop:step2n}.

\noindent Part {\rm \bf c)} follows from part {\rm \bf a)} and the observation that both $2n+1=2(3k+2)+1=6k+5$ and $2n=6k+4$ are not divisible by 3.

\noindent Part {\rm \bf d):} Type {\rm (II)} is excluded by part {\rm \bf a)}, and type {\rm (IV)} is excluded because $2n=6k+2$ is not divisible by 3. So $\Pi(n)$ is either $2n-1$ (type {\rm (I)}), or $(2n+1)/3$ (type III). Here $(2n+1)/3=(6k+3)/3=2k+1$. But $2k+1=2(k+1)-1$ was already assigned to $k+1$ if $\Pi(k+1)$ has type{\rm (I)}, so then $\Pi(n)=2n-1$ has type {\rm (I)}. Conversely, if $\Pi(k+1)$ has not type {\rm (I)}, i.e., type {\rm (III)} or type {\rm (IV)}, then $\Pi(n)$ has type {\rm (III)}. \hfill $\Box$

 \bigskip

\noindent Theorem \ref{th:types} can easily be transformed into the following statement.

\begin{theorem}  \label{thm:morphism}
Let $\sigma$ on the monoid $\{1,3,4\}^*$ be the monoid morphism given by
$$\sigma(1)=114, \; \sigma(3)=314, \;\sigma(4)=314.$$
Let $s=1141143141141143143141143141141143141141\dots$ be the  fixed point of $\sigma$ with prefix equal to $1$. Then $s_n=1$ \;iff\; $\Pi(n)$ has type {\rm (I)}, \;$s_n=3$ \;iff\; $\Pi(n)$ has type {\rm (III)}, and\; $s_n=4$ \;iff\; $\Pi(n)$ has type {\rm (IV)}.
\end{theorem}

\noindent This result has several corollaries. A first corollary is a
partial self-similarity property of the permutation $\Pi$.

\medskip

\begin{theorem}  \label{th:selfsim}
Let $E$ be the sequence of natural numbers defined by
$$E(3n+1)=9n+1,\; E(3n+2)=9n+4,\; E(3n+3)=9n+6\quad {\rm for\;} n\ge 0.$$
Then $\Pi(E(n))=3\Pi(n)-2$ for $n\ge 0$.
\end{theorem}

\proof
We split the proof in the three cases for the argument modulo 3 of $E$. Note that
$$\sigma^2(1)=114114314, \; \sigma^2(3)=314114314, \;\sigma^2(4)=314114314.$$
\noindent  From this we can read off the type of $\Pi(9n+1)$, $\Pi(9n+4)$ and $\Pi(9n+6)$, using Theorem \ref{th:types}.

\medskip

\noindent The case $\Pi(9n+1)$.

\noindent Subcase $\Pi(9n+1)$ has type {\rm (I)}.

\noindent  From Theorem \ref{th:types} we then obtain that $\Pi(3n+1)$ also has type {\rm (I)}.

 We find $\Pi(9n+1)=2(9n+1)-1=18n+1$, which is equal to $3\Pi(3n+1)-2=3(2(3n+1)-1)-2=18n+1$.

\noindent  Subcase $\Pi(9n+1)$ has type {\rm (III)}.

\noindent  From Theorem \ref{th:types} we then obtain that $\Pi(3n+1)$  has type {\rm (III)}, or  type {\rm (IV)}. However, type {\rm (IV)} does not occur, since $3n+1$ is not divisible by 3. We then find

$\Pi(9n+1)=(2(9n+1)+1)/3=6n+1$, which is equal to $3\Pi(3n+1)-2=3((2(3n+1)+1)/3)-2=6n+1$.

\medskip

\noindent The case $\Pi(9n+4)$.\; In this case both $\Pi(9n+4)$ and $\Pi(3n+2)$ have type {\rm (I)}. We then find $\Pi(9n+4)=2(9n+4)-1=18n+7$, which is equal to $3\Pi(3n+2)-2=3(2(3n+2)-1)-2=18n+7$.

\medskip

\noindent The case $\Pi(9n+6)$. \; In this case both $\Pi(9n+6)$ and $\Pi(3n+3)$ have type {\rm (IV)}. We then find
$\Pi(9n+6)=2(9n+6)/3=6n+4$, which is equal to $3\Pi(3n+3)-2=3(2(3n+3)/3)-2=6n+4$. \hfill $\Box$

\bigskip

\noindent {\bf Application}\; Consider the OEIS sequence A026186, with name 		``a(n) = (1/3)*(s(n) + 2), where s(n) is the n-th number congruent to 1 mod 3 in A026136.", with data ``$1, 3, 2, 7, 9, 4, 5, 15, 6, 19, 21, 8, 25, 27, 10, 11,\dots$".
In the 	COMMENTS  we find: ``Is this a duplicate of A026136? - R. J. Mathar, Aug 26 2019",  and ``A026186(n) = A026136(n) for n $<$= 10$\wedge$7. - Sean A. Irvine, Sep 19 2019". Translated to our setting, we have to show that\\[-.3cm]
$${\rm \bf (A)}\quad\Pi(n) \equiv 1 \mod 3\; \Leftrightarrow\; n\in\{E(k):k \ge 0\},
    \;{\rm and\;}{\rm \bf (B)}\quad  \tfrac13(\Pi(E(n)+2) = \Pi(n)\quad{\rm for}\; n\ge 0.$$\\[-.5cm]
One side of  Part {\bf (A)} follows from Theorem \ref{thm:morphism}, and the formula for $\sigma^2$ above:

$\Pi(9n+2)$ has type {\rm (I)}\quad $\Rightarrow$ \;$\Pi(9n+2)=2(9n+2)-1 = 18n +3 \equiv 0 \mod 3.$

$\Pi(9n+3)$ has type {\rm (IV)} $\Rightarrow$ \;$\Pi(9n+3)=2(9n+3)/3 = 6n +2 \equiv 2 \mod 3.$

$\Pi(9n+5)$ has type {\rm (I)}\quad $\Rightarrow$ \;$\Pi(9n+5)=2(9n+5)-1 = 18n +9 \equiv 0 \mod 3.$

$\Pi(9n+7)$ has type {\rm (III)} $\Rightarrow$ \;$\Pi(9n+7)=(2(9n+7)+1)/3 = 6n +5 \equiv 2 \mod 3.$

$\Pi(9n+8)$ has type {\rm (I)}\quad $\Rightarrow$ \;$\Pi(9n+8)=2(9n+2)-1 = 18n +15 \equiv 0 \mod 3.$

$\Pi(9n+9)$ has type {\rm (IV)} $\Rightarrow$ \;$\Pi(9n+9)=2(9n+9)/3 = 6n +6 \equiv 0 \mod 3.$

 Since all six are not equal to 1 modulo 3, but they are for the remaining three cases modulo 9 by Theorem \ref{th:selfsim}, this proves part {\bf (A)}. Part {\rm \bf (B)} follows directly from Theorem \ref{th:selfsim}.

\bigskip

Let $R_{\rm pos}=A026138$ be the sequence of positions of records in $\Pi$:
$$R_{\rm pos} = 1, 2, 4, 5, 8, 10, 11, 13, 14, 17, 20, 22, 23, 26, 28, 29, 31, 32, 35, 37, 38, 40, 41\dots$$
The sequence $R_{\rm rec}=A026139$ of records in $\Pi$ defined by $R_{\rm rec}(n)=\Pi(R_{\rm pos}(n))$ is given by
$$R_{\rm rec} = 1, 3, 7, 9, 15, 19, 21, 25, 27, 33, 39, 43, 45, 51, 55, 57, 61, 63, 69, 73, 75, 79, 81, 87, 93, 97, 99, 105, 111\dots$$

In the sequel we write $\Delta x$ for the sequence of first differences $\big (x(n+1)-x(n)\big)$ of a sequence $x$.

\medskip

\begin{proposition}\label{prop:rec} Let $\tau$ be the morphism on $\{1,2,3\}^*$ given by
$\tau(1)=12,\;\tau(2)=132,\;\tau(3)=1332.$\\
 Let $t=12132121332132\dots$ be the unique fixed point of $\tau$. \,Then $\Delta R_{\rm pos}=t$, and $\Delta R_{\rm rec}=2t$.
\end{proposition}

\proof The records in $\Pi$ are exactly given by the $\Pi(n)$ of type {\rm (I)}. So the positions of these records are given by the positions of 1 in  $s$. To obtain these, one considers the return words of the word 1 in $s$. These are $1, 14$ and $143$. Since
$$\sigma(1)=1\,14,\;  \sigma(14)=1\,143\,14,\;  \sigma(143)=1\,143\,143\,14,$$
the induced derivated morphism is equal to $\tau$, where we code the return words by their lengths.
Since we coded the return words by their lengths, this gives that $R_{\rm pos}(n+1)-R_{\rm pos}(n)=t(n)$, where $t$ is the unique fixed point of $\tau$.
The second equation follows from the fact that all records are of type {\rm (I)}, and so
$$\qquad R_{\rm rec}(n+1)-R_{\rm rec}(n)=\Pi(R_{\rm pos}(n+1))-\Pi(R_{\rm pos}(n))=2R_{\rm pos}(n+1)-1-(2R_{\rm pos}(n)-1)=2t(n).\hspace*{1cm} \Box$$

\medskip

\noindent {\bf Remark 1} The proposition states  that $A026141=\frac12 \Delta R_{\rm rec}=t$.
Also, let $A026140(n):=  \tfrac12(R_{\rm rec}(n) - 1)= 0, 1, 3, 4, 7, 9, 10, 12, 13, 16, 19, 21, 22,..,$. It is easy to see that Proposition \ref{prop:rec} gives $\Delta A026140 = t$, again.

\medskip

 \subsection{Changing the rules: $\Pi_{\rm even}$}\label{sec:even}  

 In sequence A026172 in OEIS Kimberling varies on the left-right filling procedure by adding that one always goes to the right if $n$ is even.
 The resulting permutation of ${\mathbb N}$ is denoted by $\Pi_{\rm even}$.

We still have  $L(n)=[n/2]=R(n)$, but the procedure is changed to\\[-.2cm]

\hspace*{2cm} $a(n+R(n))=n$ if $n$ even or $a(n-L(n))$ already defined, else $a(n-L(n))=n$.

\medskip

\noindent Surprisingly, changing the procedure in this way does not change the
permutation.

\begin{proposition}\label{prop:even} $\Pi_{\rm even}=\Pi$.
\end{proposition}

\proof
This follows immediately from Lemma \ref{lem:step2n}.\hfill $\Box$

\medskip

Proposition \ref{prop:even} implies equality of many pairs of sequences in OEIS,
like
   A026136 =  A026172,	
   A026137 =  A026173,	
   A026138 =  A026174,
   A026139 =  A026175,
   A026141 =  A026176,	
   A026184 =  A026208,
   A026188 =  A026212,
   A026182 =  A026206, and
   A026186 =  A026210.

\subsection{Changing the rules: $\Pi_{\rm odd}$}\label{subs:odd} 

 Still  $L(n)=[n/2]=R(n)$, but now one always goes to the right if $n$ is odd, see A026177:\\[-.3cm]

\hspace*{2cm} $a(n+R(n))=n$ if $n$ odd or $a(n-L(n))$ already defined, otherwise $a(n-L(n))=n$.\\[-.3cm]

\noindent The resulting permutation of ${\mathbb N}$ is denoted by $\Pi_{\rm odd}$. Thus
$$\Pi_{\rm odd}=1,4,2,3,10,12,5,16,6,7, 22, 8, 9, 28, 30, 11, 34, 36, 13, 40, 14, 15, 46, 48, 17, 52, 18, 19,\dots$$

\noindent  One verifies that the $\Pi_{\rm odd}(n)$ have the same four types as in Section \ref{subs:standard}, but this time it is Type {\rm (I)} that does not occur, and the records occur at Type {\rm (II)}.

\begin{theorem}  \label{th:typesodd}
Let $n,k$ be  natural numbers. Then\\
{\rm \bf a)}\; $\Pi_{\rm odd}(n)$ never has type {\rm (I)}.\\
{\rm \bf b)}\; If\, $n=3k+1$ then $\Pi_{\rm odd}(n)$ has type {\rm (III)}. \\
{\rm \bf c)}\; If\, $n=3k+2$ then $\Pi_{\rm odd}(n)$ has type {\rm (II)}.\\
{\rm \bf d)}\; If\, $n=3k$ then $\Pi_{\rm odd}(n)$ has type {\rm (II)} if $\Pi_{\rm odd}(k)$ has type {\rm (II)} and $\Pi_{\rm odd}(n)$ has type {\rm (IV)} if $\Pi_{\rm odd}(k)$ has type {\rm (III)} or  type {\rm (IV)}.
\end{theorem}

\proof
Part {\rm \bf a)}: $\Pi_{\rm odd}(n)$ is never equal to $2n-1$, since the algorithm always chooses the  right side at the odd numbers.

\noindent Part {\rm \bf b)},  i.e., $\Pi_{\rm odd}(3k+1)=2k+1$ is also forced by the `always to the right at the odd numbers' rule: $2k+1+\lfloor (2k+1)/2\rfloor=3k+1$.

\noindent Part {\rm \bf c)}, follows from part {\rm \bf a)} and the observation that both $2n+1=2(3k+2)+1=6k+5$ and $2n=6k+4$ are not divisible by 3.

\noindent Part {\rm \bf d):} Type {\rm (I)} is excluded by part {\rm \bf a)}, and type {\rm (III)} is excluded because $2n+1=6k+1$ is not divisible by 3. So $\Pi_{\rm odd}(n)$ is either $2n$ (type {\rm (II)}), or $2n/3$ (type IV). Here $2n/3=6k/3=2k$. But $2k$ would have been assigned to $k$ if $\Pi_{\rm odd}(k)$ has type{\rm (II)}, so if $\Pi_{\rm odd}(n)$ has type {\rm (IV)} then $\Pi_{\rm odd}(k)$ has type (III) or (IV). Conversely, if $\Pi_{\rm odd}(k)$  has type (III) or (IV),  then $\Pi_{\rm odd}(n)$ must have type {\rm (III)}: in $2k$ you do not go the left, since $2k-\lfloor 2k/2\rfloor=k$, and $\Pi_{\rm odd}(k)$ has already got a value. Going to the right in $2k$ yields $\Pi_{\rm odd}(2k+k)=2k$. \hfill $\Box$

 \medskip

\noindent Let $\Pi_{\rm odd}(1)=1$ have Type {\rm (III)} by definition. Then Theorem \ref{th:typesodd} can easily be transformed into the following result.

\begin{theorem}  \label{thm:morphodd}
Let $\sigma$ on the monoid $\{2,3,4\}^*$ be the monoid morphism given by
$$\sigma(2)=322, \; \sigma(3)=324, \;\sigma(4)=324.$$
Let $s_{\rm odd}=32432232432432232232432232432\dots$ be the unique fixed point of $\sigma$. Then $s_{\rm odd}(n)=2$ \;iff\; $\Pi_{\rm odd}(n)$ has type {\rm (II)}, \;$s_{\rm odd}(n)=3$ \;iff\; $\Pi_{\rm odd}(n)$ has type {\rm (III)},
and\; $s_{\rm odd}(n)=4$ \;iff\; $\Pi_{\rm odd}(n)$ has type {\rm (IV)}.
\end{theorem}

\bigskip

\noindent Let $R_{\rm opos}=A026179$ be the sequence of positions of records in $\Pi_{\rm odd}$:
$$R_{\rm opos} = 1,2, 5, 6, 8, 11, 14, 15, 17, 18, 20, 23, 24, 26, 29, 32, 33, 35, 38, 41, 42, 44\dots$$
The sequence $R_{\rm orec}=A026180$ of records in $\Pi$ defined by $R_{\rm orec}(n)=\Pi_{\rm odd}(R_{\rm opos}(n))$ is given by
$$R_{\rm orec} =1, 4, 10, 12, 16, 22, 28, 30, 34, 36, 40, 46, 48, 52, 58, 64, 66, 70, 76, 82, 84,\dots$$

\begin{proposition}\label{prop:recodd} Let $\tau$ be the morphism on $\{1,2,3\}^*$ given by
$\tau(1)=12,\;\tau(2)=312,\;\tau(3)=3312.$\\
 Let $t=12312331212312331233121231212312\dots$ be the  fixed point of $\tau$ starting with 1. Let $T$ denote the shift operator. Then $T(\Delta R_{\rm opos})=T^2(t)$, and $T(\Delta R_{\rm orec})=2\,T^2(t)$.
\end{proposition}

\proof
The records in $\Pi_{\rm odd}$ are exactly given by the $\Pi(n)$ of type {\rm (II)}. So the positions of these records are given by the positions of 2 in  $s_{\rm odd}$. To obtain these, one considers the return words of the word 2 in $s_{\rm odd}$. These are $2, 23$ and $243$. We have $$\sigma(2)=32\,2,\;  \sigma(23)=32\,23\,24\;  \sigma(243)=32\,23\,243\,24.$$
If we code the return words by their lengths, then the induced descendant morphism $\sigma_{\rm desc}$, obtained by conjugating $\sigma$ with the word 3, (see the paper \cite{Hol-Zam})  is equal to
$$ \sigma_{\rm desc}(1)= 12, \;  \sigma_{\rm desc}(2)= 123,\; \sigma_{\rm desc}(3)= 1233. $$
The `right' morphism however, is $\tau$, given by
$$\tau(1)=12, \;\tau(2)=312, \;\tau(3)=3312.$$
To prove this, avoiding confusions of symbols, we write $\tau$ on the alphabet $\{a,b,c\}$:
$$\tau(a)=ab, \;\tau(b)=cab, \;\tau(c)=ccab.$$
Let $x_\tau=abcabcc\dots$ be the fixed point of $\tau$ starting with $a$, and let $\delta$ be the morphism given by
$$\delta(a)=2,\;\delta(b)=23,\;\delta(c)=243.$$
Note that the images of $\delta$ are the return words of 2, with lengths 1,2, and 3, so the proposition will be proved if we show that
$$T^2(\delta(x_\tau))=T(s_{\rm odd})=2432232432\dots.$$
As in the paper \cite{Dekdec}, we call $\delta(x_\tau)$ a decoration of $x_\tau$. It is well-known that such a decoration is again a morphic sequence, see, e.g., the monograph by Allouche and Shallit \cite[Corollary 7.7.5]{all-shall}. We perform what is called the `natural algorithm' in the paper \cite{Dekdec}, to find the morphism and the letter-to-letter map
which yield $\delta(x_\tau)$ as a morphic sequence. Consider the alphabet $\{a, b, b', c, c', c''\}$, and define a block-substitution by
$$ a\rightarrow abb', \quad bb'\rightarrow cc'c''abb', \quad cc'c''\rightarrow cc'c''cc'c''abb'.$$
We obtain from this a morphism on $\{a, b, b', c, c', c''\}$ by splitting the images of $bb'$ and $cc'c''$ in the block-substitution (in the most efficient way):
$$ a\rightarrow abb', \quad b\rightarrow cc'c'',\quad  b'\rightarrow abb',
 \quad c\rightarrow cc'c'', \quad c'\rightarrow cc'c'',\quad c''\rightarrow abb'.$$
Here efficient means that as many as possible symbols can be merged, respecting the letter to letter map $a\rightarrow 2,\; b\rightarrow 2,\; b'\rightarrow 3,\; c\rightarrow 2,\; c' \rightarrow 4,\; c''\rightarrow 3$, obtained by identifying $a, bb'$ and $cc'c''$ with the return words 2, 23 and 243.
We merge $b$ and $c$, renaming this as $\overline{2}$, and we merge $b'$ and $c''$, renaming the resulting symbol as 3. Also renaming $a$ as 2, and $c'$ as 4, we obtain a morphism $\theta$ on the alphabet $\{2,\overline{2},3,4\}$, and a letter to letter map which is the identity, except that $\overline{2}$ is mapped to 2. Thus $\theta$ is given by
$$\theta(2)=2\overline{2}3, \quad, \theta(\overline{2})=\overline{2}43, \quad,\theta(3)=2\overline{2}3, \quad\theta(4)=\overline{2}43.$$
To obtain the claim above we have to shift the fixed point of $\theta$ by 2. The way to generate that sequence is to pass to the words of length 3 of the language of $\theta$, and to project these on the third symbol. There are 8  words of length 3:
$$1:=2\overline{2}3, \;2:=\overline{2}32,\; 3:=\overline{2}3\overline{2},
    \; 4:=\overline{2}43,\; 5:= 32\overline{2},\; 6:=3\overline{2}4, \; 7:=432,\; 8:=43\overline{2}.$$
The induced 3-block morphism (cf. the paper \cite{Dekking-JIS}) is given by
$$1 \rightarrow 136, \; 2 \rightarrow 475, \;3\rightarrow 475,
    \; 4\rightarrow 486, \;  5\rightarrow 125, 6\rightarrow 136, \; 7\rightarrow 475,\; 8\rightarrow 475.$$
The projection on the third symbol is given by
$$1\rightarrow 3,\; 2 \rightarrow 2, \;3\rightarrow \overline{2},\; 4 \rightarrow 3, \;5\rightarrow \overline{2},\; 6 \rightarrow 4,\; 7\rightarrow 2,\; 8 \rightarrow \overline{2}.$$
We see that we can consistently merge 7 and 2, and also 8 and 3. The 3-block morphism with these symbols merged is then given by
$$1 \rightarrow 136, \; 2 \rightarrow 425, \;3\rightarrow 425,
    \; 4\rightarrow 436, \;  5\rightarrow 125,\; 6\rightarrow 136.$$
We now see that we can merge 1 and 4 (both map to 3), and then also 3 and 5 (both map to $\overline{2}$), which leads to the morphism
$$1\rightarrow 136, \; 2 \rightarrow 123, \;3\rightarrow 123,\;  6\rightarrow 136.$$
Equivalently, on the 'third symbol' alphabet $\{2,\overline{2},3,4\}$:
$$3\rightarrow 3\overline{2}4, \; 2 \rightarrow 32\overline{2}, \;\overline{2}\rightarrow 32\overline{2}, \; 4\rightarrow 3\overline{2}4.$$
We can then merge 2 and $\overline{2}$ to 2, obtaining the morphism
$$3\rightarrow 324, \; 2 \rightarrow 322,  \; 4\rightarrow 324,$$
which is nothing else than the morphism $\sigma$ generating $s_{\rm odd}$. Since the first 2 in $s_{\rm odd}$ occurs at the second index, we thus proved that $T^2(\delta(x_\tau))=T(s_{\rm odd})=2432232432\dots.$

\medskip

The second equation follows from the fact that all records are of type {\rm (II)}, and so
$$\qquad R_{\rm orec}(n+1)-R_{\rm orec}(n)=\Pi_{\rm odd}(R_{\rm opos}(n+1))-\Pi_{\rm odd}(R_{\rm opos}(n))=2R_{\rm opos}(n+1)-2R_{\rm opos}(n)=2\,t(n).\hspace*{2cm} \Box$$

\subsection{Changing the rules: rule 42}\label{subs:III} 

 Here  $L(n)=[(n+1)/2]=R(n)$, and
  $a(n-L(n))=n$ if $a(L)$ not yet defined, else $a(n+R(n))=n$.

\noindent The resulting permutation of ${\mathbb N}$ is denoted by $\Pi_{\rm 42}$.
 It is sequence A026142 in the encyclopedia \cite{oeis}, given by
 $$\Pi_{\rm 42}=1, 4, 2, 8, 3, 12, 14, 5, 6, 20, 7, 24, 26, 9, 10, 32, 11, 36, 38, 13, 42, 44, 15, 16,\dots.$$

 \medskip

\noindent  One verifies that the $\Pi_{\rm 42}(n)$ have four types slightly different from those  in Section \ref{subs:standard}:\\[-.6cm] \begin{align*}
{\rm (I)}  &\quad \Pi_{\rm 42}(n)\; {\rm odd},\;\,  \; \Pi_{\rm 42}(n)>n,  \; \Pi_{\rm 42}(n) = 2n+1\\
{\rm (II)} &\quad \Pi_{\rm 42}(n)\; {\rm even},\,   \; \Pi_{\rm 42}(n)>n,  \; \Pi_{\rm 42}(n) = 2n\\
{\rm (III)} &\quad \Pi_{\rm 42}(n)\; {\rm odd},\;\, \; \Pi_{\rm 42}(n)<n,  \; \Pi_{\rm 42}(n) = (2n-1)/3\\
{\rm (IV)} &\quad \Pi_{\rm 42}(n)\; {\rm even},\,   \; \Pi_{\rm 42}(n)<n,  \; \Pi_{\rm 42}(n) = 2n/3.
\end{align*}

\noindent The next task is to determine the type of all entries of $\Pi$.

\begin{theorem}  \label{th:types-chIII}
Let $n,k$ be  natural numbers. Then\\
{\rm \bf a)}\; $\Pi_{\rm 42}(n)$ never has type {\rm (I)}.\\
{\rm \bf b)}\; If\, $n=3k+1$ then $\Pi_{\rm 42}(n)$ has type {\rm (II)}. \\
{\rm \bf c)}\; If\, $n=3k+2$ then $\Pi_{\rm 42}(n)$ has type {\rm (III)}.\\
{\rm \bf d)}\; If\, $n=3k+3$ then $\Pi_{\rm 42}(n)$ has type {\rm (II)} if $\Pi_{\rm 42}(k+1)$ has type {\rm (II)} and $\Pi_{\rm 42}(n)$ has type {\rm (IV)} if $\Pi_{\rm 42}(k+1)$ has type {\rm (III)} or  type {\rm (IV)}.
\end{theorem}

\proof
Part {\rm \bf a)}: $\Pi_{\rm 42}(n)$ is never equal to $2n$, since either $n$ has been assigned the value $2n+1$ at step $2n+1$:
$2n+1-\lfloor (2n+2)/2\rfloor=n$, or $n$ has been assigned a smaller value at an earlier step.

\noindent Part {\rm \bf b)}  follows from part {\rm \bf a)} and the observation that both $2n-1=2(3k+1)-1=6k+1$ and $2n=6k+2$ are not divisible by 3.

\noindent Part {\rm \bf c)} is a rephrasing of the version of Proposition \ref{prop:step2n} for rule 42: $\Pi(3n+2)=2n+1$, which is a consequence of the version of Lemma \ref{lem:step2n} for rule 42 (`if the procedure is at step $2n+1$, then the numbers $1,\dots,n$ have already been assigned a $\Pi$-value.')

\noindent Part {\rm \bf d):} Type {\rm (I)} is excluded by part {\rm \bf a)}, and type {\rm (III)} is excluded because $2n-1=6k+5$ is not divisible by 3. So  $\Pi_{\rm 42}(n)$ is either $2n$ (type {\rm (II)}), or $2n/3$ (type IV). Here $2n/3=(6k+6)/3=2k+2$. But $2k+2=2(k+1)$ was already assigned to $k+1$ if $\Pi(k+1)$ has type {\rm (II)}, so   $\Pi_{\rm 42}(n)=2n/3$ has type {\rm (II)} when $k+1$ has $\Pi(k+1)$ has type{\rm (II)}. Conversely, if $\Pi(k+1)$ has not type {\rm (II)}, i.e., type {\rm (III)} or type {\rm (IV)}, then $\Pi_{\rm 42}(n)$ has type {\rm (IV)}. \hfill $\Box$

 \medskip

 For the next result we need a fifth type. We say $\Pi_{\rm 42}(n)$ has type (V) if $\Pi_{\rm 42}(n) = n$. Actually the only $\Pi_{\rm 42}(n)$ of type (V) is $\Pi_{\rm 42}(1)$.
This way one arrives at the following result.

\begin{theorem}  \label{thm:morphismodd}
Let $\sigma$ on the monoid $\{2,3,4,5\}^*$ be the monoid morphism given by
$$\sigma(2)=232, \; \sigma(3)=234, \;\sigma(4)=234, \;\sigma(5)=524.$$
Let $s=5, 2, 4, 2, 3, 2, 2, 3, 4, 2, 3, 2, 2, 3, 4, 2, 3, 2, 2, 3, 2, 2, 3, 4, 2, 3, 4, 2, 3, 2, 2, 3, 4\dots$  be the fixed point of $\sigma$ starting with $5$. Then for all $n>1$: $s_n=2$ \;iff\; $\Pi_{\rm 42}(n)$ has type {\rm (II)}, \;$s_n=3$ \;iff\; $\Pi_{\rm 42}(n)$ has type {\rm (III)}, and\; $s_n=4$ \;iff\; $\Pi_{\rm 43}(n)$ has type {\rm (IV)}.
\end{theorem}

\bigskip

\noindent Let $R_{\rm pos}=A026144$ be the sequence of positions of records in $\Pi_{\rm 42}$:
$$R_{\rm pos} = 1, 2, 4, 6, 7, 10, 12, 13, 16, 18, 19, 21, 22, 25, 28, 30, 31, 34, 36, 37, 39\dots$$
The sequence $R_{\rm rec}=A026145$ of records in $\Pi$ defined by $R_{\rm rec}(n)=\Pi_{\rm 42}(R_{\rm pos}(n))$ is given by
$$R_{\rm rec} =1, 4, 8, 12, 14, 20, 24, 26, 32, 36, 38, 42, 44, 50, 56, 60, 62, 68, 72, 74, 78, 80,\dots$$
Except for the first one, the records are always even, because they are generated by type (II) elements of the permutation.

\begin{proposition}\label{prop:rec-chIII} Let $\tau$ be the morphism on $\{1,2,3,4\}^*$  given by
$$\tau(1)=21,\;\tau(2)=213,\;\tau(3)=2133,\;\tau(4)=4213.$$
 Let $t=421321321213321321213321321213212\dots$ be the  fixed point of $\tau$ starting with $4$. \,Then $\Delta R_{\rm pos}(n+1)=t(n)$,  for $n \ge 2$ and $\Delta R_{\rm rec}(n+1)=2t(n)$,  for $n\ge 2$.
\end{proposition}

\proof The records in $\Pi$ are exactly given by the $\Pi(n)$ of type {\rm (II)}. So the positions of these records are given by the positions of 2 in  $s$. To obtain these, one considers the return words of the word 2 in $s$. These are $2, 23$, $24$ (only at the beginning) and $234$. Since
$$\sigma(2)=23\,2,\;  \sigma(23)=23\,2\,234,\;  \sigma(24)=23\,2\,234,\; \sigma(234)=23\,2\,234\,234,$$
the induced derivated morphism is equal to $\tau$, where we code the return words by their lengths, and we added a `starting' letter $4$ at the beginning.
Since we coded the return words by their lengths, this gives that $R_{\rm pos}(n+1)-R_{\rm pos}(n)=t(n)$, for $n\ge 2$ where $t$ is the  fixed point of $\tau$ starting with $4$.
The second equation follows from the fact that all records are of type {\rm (II)}, and so
$$\qquad R_{\rm rec}(n+1)-R_{\rm rec}(n)=\Pi(R_{\rm pos}(n+1))-\Pi(R_{\rm pos}(n))=2R_{\rm pos}(n+1)-2R_{\rm pos}(n)=2t(n).\hspace*{1cm} \Box$$

\subsection{Comparing two permutations}\label{subs:twoper}

 The  permutation $\Pi_{\rm 36}:=\Pi$ from Section \ref{subs:standard} and $\Pi_{\rm 42}$ from Section \ref{subs:III} have many entries in common, as observed in OEIS sequence A026222 = $1, 3, 9, 15, 24, 27, 33, 42, 45, 51, 60, 69, 72, \dots$, which gives the first 66 numbers $n$ satisfying $\Pi_{\rm 36}(n)=\Pi_{\rm 42}(n)$. Here we prove that there are infinitely many of such entries. More precisely, let  this `coincidence' sequence  be $C$ given by: $n$ is in $\{C(k): k\in \mathbb{N}\}$  if and only if $\Pi_{\rm 36}(n)=\Pi_{\rm 42}(n)$. We will show that the difference sequence $\Delta C$ is 3-automatic.

\begin{theorem}  \label{thm:coinc}
Let $C$  be the `coincidence' sequence of $\Pi_{\rm 36}$ and $\Pi_{\rm 42}$. Then $C=I_4$, where $I_4$ is the sequence of natural numbers with type {\rm (IV)} in $\Pi_{\rm 42}$.
\end{theorem}

\proof  To generate  the `coincidence' sequence  $C$, we have to consider the product of the two morphisms generating the types of  the two permutations $\Pi_{\rm 36}$ and $\Pi_{\rm 42}$. This product is defined on the set of product symbols $\{1,2,3,4\}\times\{1,2,3,4\}$. Not all product symbols occur. By writing out the images of the product symbols under the product substitutions one produces a list of the relevant ones:\\[-.1cm]
$$\PQ{1}{5},\; \PQ{1}{2},\;\PQ{4}{4},\; \PQ{1}{3},\; \PQ{4}{2},\; \PQ{3}{2}.  $$
One has, for example,\\[-.1cm]
$$\PQ{1}{5}\rightarrow \PQ{1}{5} \PQ{1}{2} \PQ{4}{4},\;\PQ{1}{2}\rightarrow \PQ{1}{2} \PQ{1}{3} \PQ{4}{2},
   \; \PQ{4}{4}\rightarrow\PQ{3}{2}\PQ{1}{3}\PQ{4}{4}.$$
 From the observation that a coincidence $\pq{4}{4}$ occurs if and only if a $4$ occurs in the `lower' sequence $\Pi_{\rm 42}$, the theorem now follows. \hfill $\Box$

\medskip

The next result identifies the sequence $I_4$ in terms of its first differences. We do need again an extra symbol $4$ to deal with the beginning of $I_4$.

\begin{proposition}\label{prop:coinc} Let $\kappa$ be the morphism on $\{1,2,3,4\}^*$  given by
$$\kappa(1)=12,\;\kappa(2)=123,\;\kappa(3)=1233,\;\kappa(4)=423.$$
 Let $k=4231231233\dots$ be the  fixed point of $\kappa$ starting with $4$. \,Then $\Delta I_4(n+1)=\lambda(k(n))$ for $n\ge 1$,  where $\lambda$ is the letter-to-letter map $1\rightarrow 3, 2\rightarrow 6, 3 \rightarrow 9, 4\rightarrow 6$.
\end{proposition}

\proof
Let $s_{\rm 42}=5, 2, 4, 2, 3, 2, 2, 3, 4, 2, 3, 2, 2, 3, 4, 2, 3, 2, 2, 3,$  be the sequence of types of $\Pi_{\rm 42}$ given by Theorem \ref{thm:morphismodd}, which gives $s_{\rm 42}$ as fixed point of the morphism  $\sigma$ defined by
$$\sigma(2)=232, \; \sigma(3)=234, \;\sigma(4)=234, \;\sigma(5)=524.$$
The  three return words of the word 4 in  $s$ are   $a:=423,\,b:=423223$ and $c:=423223223$. We  can not apply the descendant algorithm of the paper \cite{Hol-Zam}, since the sequence $s$  is not a minimal sequence, not only because of the unique appearance of $s_{\rm 42}(1)=5$ at the beginning, but also because of the unique
appearance of $s_{\rm 42}(2)s_{\rm 42}(3)=24$, $s_{\rm 42}(3)\dots s_{\rm 42}(14)=bb$, etc..

We compute the images of the return words under $\sigma$.
$$\sigma(a)=23\,b\,4,\;  \sigma(b)=23\,bc\,4,\;\sigma(c)=23\,bcc\,4.$$
We see from this that if we write
$$s_{\rm 42}=\sigma(s_{\rm 42})=52\,b'bcabcabcc\dots,$$
where $b'=423223$, then the sequence $s_{\rm 42}(3)s_{\rm 42}(4)s_{\rm 42}(5)\dots$ will be a fixed point
of the morphism given by
$$a\rightarrow ab, \; b\rightarrow abc,  \; b'\rightarrow b'bc, \; c\rightarrow abcc.$$
The proposition now follows by changing the alphabet $\{a,b,c,b'\}$ to  $\{1,2,3,4\}$,
and noting that $|a|=3, |b|=|b'|=6, |c|=9$. \hfill $\Box$

\medskip

We remark that the main, primitive part of $\kappa$ is given by A106036.

\subsection{Changing the rules: $\Pi(1)$}\label{subs:twoper}

One might get rid of the requirement $\Pi(1):=1$ by defining a new permutation $\Pi_\oplus$ by
$$\Pi_\oplus(n):=\Pi(n+1)-1\quad {\rm for\;} n=1,2,\dots.$$
For the rule of Section \ref{subs:standard} given by  $L(n)=\lfloor n/2\rfloor, \;R(n) = \lfloor n/2\rfloor,$ this transforms
$$\Pi=1,3,2,7,9,4,5,15,6,19,21,8,25,27,10,11,33,12,13,39,14,43,45,16,17,51,\dots$$
into
$$\Pi_\oplus=2, 1, 6, 8, 3, 4, 14, 5, 18, 20, 7, 24, 26, 9, 10, 32, 11, 12, 38, 13, 42, 44, 15, 16, 50\dots.$$
There is a remarkable connection with the permutation $\Pi_{\rm odd}$. In the following proposition the sequence $s$ appears as A026215 in OEIS.

\begin{proposition}\label{prop:oplus} Let $\Pi=\Pi_{\rm 36}$ be the permutation of Section \ref{subs:standard}, and let $\Pi_{\rm odd}$ be the permutation of Section \ref{subs:odd}, then $\Pi_\oplus(n)$ is the position of $n$ in the sequence $(s(n)/2)$, where $s(n)$ is the $n^{\rm th}$ even number in $\Pi_{\rm odd}$.
\end{proposition}

\proof
The positions of the even numbers in $\Pi_{\rm odd}$ are given by the entries 2 and 4 of the fixed point $s_{\rm odd}=32432232432432232232432232432\dots$ of the  morphism $\sigma_{\rm odd}$ given by
$$\sigma_{\rm odd}(2)=322, \; \sigma_{\rm odd}(3)=324, \;\sigma_{\rm odd}(4)=324.$$
See Theorem \ref{thm:morphodd}. We see from the form of $\sigma_{\rm odd}$ that the positions of the even numbers in $\Pi_{\rm odd}$ are then exactly given by the union of the  two arithmetic sequences $(3n+2:n\ge 0)$ and $(3n+3:n\ge 0)$. This implies that
$$s(2n+1) = \Pi_{\rm odd}(3n+2)\quad {\rm and\quad}s(2n+2) = \Pi_{\rm odd}(3n+3)  \;{\rm for\quad}n=0,1,2,\dots .$$
So we have to prove for $j=1,2$ that
$$ \Pi_\oplus^{-1}(2n+j)=\tfrac12\Pi_{\rm odd}(3n+j+1)
\;\Leftrightarrow\; \Pi_\oplus\big(\tfrac12\Pi_{\rm odd}(3n+j+1)\big)=2n+j.$$
Here the case $j=1$ is the simplest: we see from the form of $\sigma_{\rm odd}$ that all $n$ which are 2\, modulo\, 3 have type (II), hence $\Pi_{\rm odd}(3n+2)/2=(6n+4)/2=3n+2$, and also $\Pi_\oplus(3n+2)=\Pi(3n+3)-1=2n+1$, since all multiples of 3 in the permutation $\Pi$ have type (IV), by Theorem \ref{th:types}.

The case $j=2$ is similar, but more involved. We must see that  $\Pi_\oplus(\Pi_{\rm odd}(3n+3)/2)=2n+2$, which holds, replacing $n+1$ by $n$ in this equation, if and only if $$\Pi\big(\tfrac12\Pi_{\rm odd}(3n)+1\big)=2n+1.$$
There are two possible types for $\Pi_{\rm odd}(3n)$: (A) type (IV), and (B) type (II).

In case (A) we have $\Pi_{\rm odd}(3n)=2n$, so then we have to see that $\Pi(n+1)=2n+1$, which means that $\Pi(n+1)$ should have type (I). Using Theorem \ref{th:types}, part d), this holds if and only if $\Pi(3n+1)$ has type (I).

In case (B) we have $\Pi_{\rm odd}(3n)=6n$, so then we have to see that $\Pi(3n+1)=2n+1$, which means that $\Pi(3n+1)$ should have type (III). We should therefore prove the following.
\begin{align*}
{\rm (A)} & \quad\Pi(3n+1) \;{\rm has\; Type\; (I)}\quad {\rm iff}\quad \Pi_{\rm odd}(3n) \;{\rm has\; Type\; (IV)}\;\\
{\rm (B)} & \quad\Pi(3n+1) \;{\rm has\; Type\; (III)}\;\;   {\rm iff}\quad \Pi_{\rm odd}(3n)\;{\rm has\; Type\; (II)}\;.
\end{align*}
Let $\sigma=:\sigma_{\rm 36}$ be the morphism generating the sequence $s_{\rm 36}$ of types of $\Pi$. The sequence $Ts_{\rm 36}=s_{\rm 36}(2)s_{\rm 36}(3)\dots$ is also 3-automatic, and a simple computation yields that $Ts_{\rm 36}$ is the fixed point of the morphism $\tau$ on $\{1,3,4\}$ given by
$$\tau(1)=141, \tau(3)=143, \tau(4)=143.$$
Now consider the product morphism $\tau\times \sigma_{\rm odd}$ on the relevant subset of $\{1,3,4\}\times\{2,3,4\}$.
Starting with the first symbol $(1,3)$ of the  fixed point of the product morphism, one finds quickly that these relevant symbols are
$(1,3), \,(4,2),\, (1,4)$ and $(3,2)$. It is also easy to see that at entries which are a multiple of 3, only the two symbols
$(1,4)$ and $(3,2)$ occur. This proves (A) and (B). \hfill $\Box$

\section{Conclusion}

 We have introduced a frame work to analyse permutations generated by a left-right filling algorithm. This has lead to a compact description of sequences A026136, A026142, A026172, A026177 in OEIS. We leave the permutation given in A026166 to the interested reader. There are numerous relations that can be derived from this description as we showed, for example, at the end of Section \ref{sec:even}. Here is one conjecture on the basic permutation $\Pi=\Pi_{\rm 36}$ in Section \ref{subs:standard}:
   $R_{\rm pos}$ and $R_{\rm rec}$ seem to be disjoint sequences. The complement of their union: $\{6,12,16,18,24,30,\ldots\}$, divided by 2:   $\{3,6,8,9,12,15,\ldots\}$ is equal to  A189637, the  positions of 1 in A116178, where A116178 is	Stewart's choral sequence,  unique fixed point of the morphism $0\rightarrow001, 1\rightarrow011$.

%
%
%
%
%
%
%


\end{document}